\documentclass[12pt, reqno]{amsart}
\usepackage{amsmath, amsthm, amscd, amsfonts, amssymb, graphicx, color}
\usepackage[bookmarksnumbered, colorlinks, plainpages]{hyperref}

\newtheorem{theorem}{Theorem}[section]
\newtheorem{lemma}[theorem]{Lemma}

\newtheorem{corollary}[theorem]{Corollary}
\theoremstyle{definition}
\newtheorem{definition}[theorem]{Definition}
\newtheorem{example}[theorem]{Example}
\newtheorem{Notation}[theorem]{Notation}
\theoremstyle{remark}

\numberwithin{equation}{section}

\begin{document}

\title[Common fixed point theorems for two mappings in $b$-metric-like spaces]
{Common fixed point theorems for two mappings in $b$-metric-like
spaces}

\author[B. Mohebbi Najmabadi and T.L. Shateri ]{ B. Mohebbi Najmabadi and T.L. Shateri }
\address{Behroz Mohebbi Najmabadi \\Department of Mathematics and Computer
Sciences, Hakim Sabzevari University, Sabzevar.IRAN }
\email{\rm behrozmohebbi1351@gmail.com; behruz\_mohebbi@yahoo.com}
\address{Tayebe Lal Shateri \\ Department of Mathematics and Computer
Sciences, Hakim Sabzevari University, Sabzevar, P.O. Box 397, IRAN}
\email{ \rm t.shateri@gmail.com; t.shateri@hsu.ac.ir}

\thanks{*The corresponding author:
t.shateri@gmail.com; t.shateri@hsu.ac.ir  (Tayebe Lal Shateri)}
 \maketitle

\begin{abstract}
In this paper, we prove some common fixed point theorems for two
self-mappings in  $b$-metric-like spaces. \vskip 3mm
\end{abstract}
\section{\bf Introduction}\vskip 2mm
 Fixed point theory is an important and actual topic of nonlinear
analysis.  During the last four decades fixed point theorem has
undergone various generalizations either by relaxing the condition
on contractivity or withdrawing the requirement of completeness or
sometimes even both, but a very interesting generalization was
obtained by changing the structure of the space according to this
argument. Matthews \cite{SM} introduced the notion of partial metric
space. The notions of metric-like spaces \cite{AA} and $b$-metric
spaces \cite{AY,BA,BO} were introduced, which are generalizations of
metric spaces. After that fixed point results in these spaces have
been studied by many authors \cite{AAS,AN,EG1,EG2,EG3,SC,RA,ZO1,ZO2}. Recently, Alghamdi et
al. \cite{MA} introduced $b$-metric- like spaces and gave some fixed
point results in such spaces. \\In this paper, we prove some common
fixed point theorems for two self-mappings in \\$b$-metric-like spaces.
In order to do this, we present the necessary definitions and
results in $b$-metric-like Spaces, which will be useful for the rest
of the paper. For more details, we refer to \cite{MA,AAS}.
\begin{definition}\label{D}
 A mapping $P:X\times X\to \Bbb R$ is said to be a partial metric
on $X$ if for any $x,y,z\in X$, the following conditions hold:\\
$(P_1)\;x=y$  if and only if $P(x,x)=P(y,y)=P(x,y)$ \\
$(P_2)\;P(x,x)\leq P(x,y)$ \\
$(P_3)\;P(x,y)=P(y,x)$ \\
$(P_4)\;P(x,z)\leq P(x,y)+P(y,z)-P(y,y)$.\\
Then the pair $(X,P)$ is called a partial metric space.
\end{definition}
\begin{definition}
 A mapping $\sigma:X\times X\to \Bbb R$ is said to be a metric-like
on $X$ if for any $x,y,z\in X$, the following conditions hold:\\
$(\sigma_1)\; \sigma(x,y)=0\Rightarrow x=y$\\
$(\sigma_2)\; \sigma(x,y)=\sigma(y,x)$\\
$(\sigma_3)\; \sigma(x,z)\leq\sigma(x,y)+\sigma(y,z)$.\\
Then the pair $(X,\sigma)$ is called a  metric-like space.
\end{definition}
\begin{example}
Let $X=\{0,1\}$ and
\begin{equation*}
\sigma(x,x)=\left\{
\begin{array}{cc}
2 & $ x=y=0$ \\
1 &  $ otherwise $
\end{array} \right.
\end{equation*}
\end{example}
The concept of $b$-metric space was introduced by Czerwik \cite{SC}.
\begin{definition}
A $b$-metric on $X$ is a mapping $D:X\times X\to [0,\infty)$ such that
for all $x,y,z\in  X$  and  constant $K\geq 1$ the following
conditions hold:\\
$(D_1)$  $D(x,y)=0\Leftrightarrow x=y$\\
$(D_2)$  $D(x,y)=D(y,x)$\\
$(D_3)$  $D(x,y)\leq K\Big[D(x,z)+D(z,y)\Big]$.\\
Then the pair $(X,D)$ is called a $b$-metric-space.
\end{definition}
\begin{definition}
A $b$-metric-like on $X$ is a function $D:X\times X\to [0,\infty)$
such that for all $x,y,z\in X$  and an constant $K\geq 1$ the
following
conditions hold:\\
$(D_1)\;  D(x,y)=0\Rightarrow x=y;$\\
$(D_2)\;  D(x,y)=D(y,x);$\\
$(D_3)\;  D(x,y)\leq K[D(x,z)+D(z,y)].$\\
Then the pair $(X,D)$ is called a $b$-metric-like space.
\end{definition}
\begin{example}
Let $X=[0,\infty)$. Define the function $D:X^2\rightarrow
[0,\infty)$ by $D(x,y)=(\sqrt{x}+\sqrt{y})^2$. Then $(X,D,2)$ is a
$b$-metric-like space with the constant $K=2$.
\end{example}
\begin{definition}
Let $(X,D)$ be a $b$-metric-like space and let $\{x_n\}$ be a sequence
 of $X$ and $x\in X$. Then $\{x_n\}$ is said to be convergent to $x$ and
 denote it by $x_n\to x$, if $\lim_{n\to\infty} D(x,x_n)=D(x,x)$.
\end{definition}
Before we state and prove the main result, we need the following
lemmas from \cite{MA}.
\begin{lemma}\cite[Lemma 2.13]{MA}
 Let $(X,D,K)$ be a $b$-metric-like space, and let
$\{x_k\}_{k=0}^n\subseteq X$. Then
\begin{eqnarray*}
D(x_n,x_0)\leq K
D(x_0,x_1)+...+K^{n-1}D(x_{n-2},x_{n-1})+K^{n-1}D(x_{n-1},x_n).
\end{eqnarray*}
\end{lemma}
\begin{lemma}\cite[Lemma 2.14]{MA}\label{lem1}
 Let $\{y_n\}$ be a sequence in a $b$-metric-like space $(X,D,K)$ such that
 $D(y_n,y_{n+1})\leq\lambda D(y_{n-1},y_n)$, for some $\lambda$ with $0<\lambda <K^{-1}$ and each $n\in \mathbb{N}$. Then
$\lim_{n,m\to \infty} D(y_m,y_n)=0.$
\end{lemma}
\begin{lemma}\cite[Proposition 2.10]{MA}\label{lem2}
 Let $(X,D,K)$be a $b$-metric-like space, and let $\{x_n\}$  be a
sequence in $ X$ such that $ \lim_{n\to\infty}D(x_n,x)=0$. Then\\
$(i)$ $x$ is unique \\
$(ii)$ $k^{-1}D(x,y) \leq \lim_{n\to\infty}D(x_n,y)\leq KD(x,y)$,
for all $y \in X$.
\end{lemma}
\begin{Notation}
Let $(X,D,K)$ be a $b$-metric-like space. Define $D^s:X^2\rightarrow
[0,\infty]$ by $D^s(x,y)=|2D(x,y)-D(x,x)-D(y,y)|$. Clearly,
$D^s(x,x)=0 $, for all $x\in X$.
\end{Notation}
\section{Main results}
In this section we prove some common fixed point results for
mappings on a $b$-metric-like space. In fact by using some ideas of
\cite{MA} we generalize fixed point results  for two maps in
$b$-metric like spaces.

\begin{theorem}\label{th1}
Let $(X,D,K)$ be a complete  $b$-metric-like space. Assume that
 $S,T:X\to X$ are onto  mappings such that
\begin{equation}\label{1}
D(Tx,Sy)\geq \Big[R+L\min\{D^s(x,Tx),D^s(y,Sy),
D^s(x,Sy),D^s(y,Tx)\}\Big]D(x,y),
\end{equation}
for all $x,y\in X$, where $R>K$ and $L\geq 0$. Then $T$ and $S$ have
a unique common fixed point.
\end{theorem}
\begin{proof}
First we show that if $T$ and $S$ have a common fixed point, and
then  the fixed point is unique. Let $x,y$ be two common fixed points
that is $Tx=Sx=x$ and $Ty=Sy=y$. If $ x\neq y $ then from (2.1) we
have
\begin{eqnarray*}
D(x,y)=D(Tx,Sy))\geq RD(x.y)>KD(x,y)\geq D(x,y)
\end{eqnarray*}
which is a contradiction. Thus $x=y$.

 Fix $x_0\in X$. Since $T,S$ are onto, so there exists $x_1\in
X$ such that $x_0=Tx_1$ and there exists a $x_2\in X$ such that
$Sx_2=x_1$. By continuing this process, we get $x_{2n}=Tx_{2n+1}$,
$x_{2n+1}=Sx_{2n+2}$, for all  $n\in \Bbb{N}\cup \{0\}$. In case
 $x_{n_0}=x_{n_0+1}$ for some $n_0\in \Bbb{N}\cup \{0\}$.
 Then it is clear that $x_{n_0}$ is a fixed point of $T$ or $S$. Now
 assume that $x_n\neq x_{n+1}$, for all $n\in\Bbb N$. From (\ref{1}) we have
\begin{align*}
D(x_{2n+1},x_{2n})&=D(Sx_{2n+2},Tx_{2n+1})\\
&\geq\Big[R+L\min\{D^s(x_{2n+1},Tx_{2n+1}),D^s(x_{2n+2},Sx_{2n+2})\\
&,
D^s(x_{2n+1},Sx_{2n+2}),D^s(x_{2n+2},Tx_{2n+1})\}\Big]D(x_{2n+1},x_{2n+2}).
\end{align*}
Which implies that $D(x_{2n+1},x_{2n})\geq RD(x_{2n+1},x_{2n+2})$,
and similarly, we get $$D(x_{2n+2},x_{2n+1})\geq RD(x_{2n},x_{2n+1}),$$
and so $ D(x_n,x_{n+1})\leq R^{-1} D(x_{n-1},x_{n})$ for all $n$.
Lemma (\ref{lem1}) implies that $$\lim_{n,m\to\infty}D(x_n,x_m)=0.$$
So $\{x_n\}$ is a Cauchy sequence. Since $(X,D,K)$ is a complete
$b$-metric-like space the sequence $\{x_n\}$ is converges to any $z$
in $X$,  and $\{x_{2n}\}, \{x_{2n+1}\}$ are also converges to $z\in
X$. Thus
\begin{align*}
\lim_{n\to\infty}D(x_n,
z)=\lim_{n\to\infty}D(x_{2n},z)=\lim_{n\to\infty}D(x_{2n+1},
z)=\lim_{n\to\infty}D(x_n,x_m)=0.
\end{align*}
Since $T$ and $ S$ are onto, there exists a $w_1\in X$ and a $w_2\in
X$ such that $Tw_1=z$ and $Sw_2=z$. From (\ref{1}) we have
\begin{align*}
D(x_{2n},z)&=D(Tx_{2n+1},Sw_2)\\&\geq\Big[R+L\min\{D^s(x_{2n+1},Tx_{2n+1}), D^s(w_2, Sw_2),D^s(x_{2n+1},Sw_2)\\&,
D^s(w_2, Tx_{2n+1})\}\Big]D(x_{2n+1},w_2).
\end{align*}
Taking limit as $n\to \infty$ in the above inequality, we get
\begin{align*}
0=\lim_{n\to\infty}D(x_{2n}, z)\geq
R\lim_{n\to\infty}D(x_{2n+1},w_2),
\end{align*}
which implies that  $ \lim_{n\to\infty}D(x_{2n+1},w_2)=0$. Then
lemma (1.10) implies that $w_2=z$ or $Sz=z$. Similarly, we have
\begin{align*}
D(x_{2n+1},z)&=D(Sx_{2n+2},Tw_1)\\&\geq\Big[R+L\min\{D^s(x_{2n+2},Sx_{2n+2}),D^s(w_1,Tw_1),D^s(w_1,Sx_{2n+2})\\&,
D^s(x_{2n+2},Tw_1)\}\Big]D(x_{2n+2},w_1).
\end{align*}
Taking limit as $n\to \infty$ in the above inequality , we get
\begin{equation}
0=\lim_{n\to\infty}D(x_{2n+1}, z)\geq
R\lim_{n\to\infty}D(x_{2n+2},w_1)
\end{equation}
which implies that $ \lim_{n\to\infty}D(x_{2n+2},w_1)=0$ then by
lemma (1.10) $w_1=z$ that is $Tz=z$.
\end{proof}
If take $L=0$ in Theorem \ref{th1} , then we get the following
corollary.
\begin{corollary}
Let $(X,D,K)$ be a complete  $b$-metric-like space. Assume that
$S,T:X\to X$ are onto  mappings, such that; $D(Tx,Sy)\geq RD(x,y)$
for all $x,y\in X$, where $R>K$. Then $T$ and $S$ have a unique
common fixed point.
\end{corollary}
\begin{example}
Let $X=[0,\infty)$. Define the function $D:X^2\rightarrow
[0,\infty)$ by $D(x,y)=(\sqrt{x}+\sqrt{y})^2$. Then $(X,D,2)$ is a
$b$-metric-like space. Define $ S,T:X\to X $ such that $ Sy=y $ and $
Tx=9x $ and $ R=3 $ Then if we get $x\geq y$
\begin{align*}
D(Tx,Sy)=(\sqrt{9x}+\sqrt{y})^{2}=9x+y+6\sqrt{xy}\geq
3x+3y+6\sqrt{xy}=3(\sqrt{x}+\sqrt{y})^2
\end{align*}
Hence $T$ and  $S$ have a unique common fixed point.
\end{example}
In the following, we suppose $\varphi:(0,\infty)\to(L^2,\infty)$ is
a function which satisfies  the condition $\varphi(t_n)\to
(L^2)^+\Rightarrow t_n\rightarrow 0$, where $L>1$. An example of
this function is $ \varphi(t)=4(1+t) $ for $ L=2 $ .
\begin{theorem}
Let $(X,D,K)$ be a complete  $b$-metric-like space, Assume that
$S,T:X\to X$ are onto  mappings, such that
\begin{equation}\label{eq4.1}
D(Tx,Sy)\geq \varPhi\Big(D(x,y)\Big)D(x,y)
\end{equation}
 for all $x,y\in X$. Then $T$ and $S$ have a unique common fixed point.
\end{theorem}
\begin{proof}
It is clear that if $ T $ and $ S $ have a  common fixed point, then
it is unique. Suppose to the contrary that $x_0\in X$. Since $T$ and $S$ are onto,
there exists $x_1\in X$ such that $x_0=Tx_1$ and there exists
$x_2\in X$ such that $Sx_2=x_1$. By continuing this process, we get
$x_{2n}=Tx_{2n+1}$ ,$x_{2n+1}=Sx_{2n+2}$ for all  $n\in \Bbb{N}\cup
\{0\}$. In case
 $x_{n_0}=x_{n_0+1}$ for some $n_0\in \Bbb{N}\cup \{0\}$.
Then it is clear that $x_{n_0}$ is a fixed point of $T$ or $S$. Now
 assume that $x_n\neq x_{n+1}$ for all $n$ from \ref{eq4.1} we get
\begin{align*}
D(x_{2n},x_{2n+1})=D(Tx_{2n+1},Sx_{2n+2})&\geq
\varphi\Big(D(x_{2n+1},x_{2n+2})\Big)D(x_{2n+1},x_{2n+2})\\&\geq
K^2D(x_{2n+1},x_{2n+2})\geq D(x_{2n+1},x_{2n+2}),
\end{align*}
and also
\begin{align*}
D(x_{2n+1},x_{2n+2})=D(Sx_{2n+2},Tx_{2n+3})&\geq
\varphi\Big(D(x_{2n+3},x_{2n+2})\Big)D(x_{2n+3},x_{2n+2})\\&\geq
K^2D(x_{2n+3},x_{2n+2})\geq D(x_{2n+3},x_{2n+2}).
\end{align*}
Then the sequence $D(x_{n},x_{n+1})$ is a decreasing sequence in
$\Bbb R^+$ and so there exists $s\geq 0$ such that
$\lim_{n\to\infty}D(x_n,x_{n+1})=s$. Now we prove  $s=0$. Suppose to
the  Contrary that $s>0$. By \ref{eq4.1} we
 deduce
 \begin{align*}
K^2\frac{D(x_{2n},x_{2n+1})}{D(x_{2n+1},x_{2n+2})}\geq
\frac{D(x_{2n},x_{2n+1})}{D(x_{2n+1},x_{2n+2})}\geq
\varphi\Big(D(x_{2n+1},x_{2n+2})\Big)\geq K^2
\end{align*}
By taking limit as $n\to \infty$ in the above inequality, we get
\begin{align*}
\lim_{n\to\infty}B\Big(D(x_{2n+1},x_{2n+2})\Big)=k^2.\\
\end{align*}
Hence, $s=\lim_{n\to\infty}D(x_{2n+1},x_{2n+2})=0$, and also
$s=\lim_{n\to\infty}D(x_{2n+3},x_{2n+2})=0$ which is a
contradiction. Hence $s=0$. \\
We shall show that $\limsup_{n,m\to\infty}D(x_n,x_m)=0$.\\
Suppose to the contrary that $ \limsup_{n,m\to\infty}D(x_n,x_m)>0$.
Thus we have
\begin{align*}
D(x_{2n},x_{2m+1})=D(Tx_{2n+1},Sx_{2m+2})\geq
\varphi\Big(D(x_{2n+1},x_{2m+2})\Big)D(x_{2n+1},x_{2m+2})\\
\end{align*}
that is
\begin{align*}
\frac{D(x_{2n},x_{2m+1})}{\varphi\Big(D(x_{2n+1},x_{2m+2})\Big)}
\geq D(x_{2n+1},x_{2m+2}).
\end{align*}
then by $(D_3)$ we get,
\begin{eqnarray*}
D(x_{2n},x_{2m+1})&\leq&
KD(x_{2n},x_{2n+1})+K^2D(x_{2n+1},x_{2m+2})+K^2D(x_{2m+2},x_{2m+1})\\
&\leq&
KD(x_{2n},x_{2n+1})+K^2\frac{D(x_{2n},x_{2m+1})}{\varphi\Big(D(x_{2n+1},x_{2m})\Big)}+K^2D(x_{2m+2},x_{2m+1}).
\end{eqnarray*}
Therefore
\begin{align*}
D(x_{2n},x_{2m+1})\leq
\Big(1-\frac{K^2}{\varphi\Big(D(x_{2n+1},x_{2m+2})\Big)}\Big)^{-1}
\Big(KD(x_{2n},x_{2n+1})+K^2D(x_{2m+2},x_{2m+1})\Big)
\end{align*}
By taking limit as $n,m\to \infty$ in the above inequality, since
$\limsup_{n,m\to\infty}D(x_{2n},x_{2m+1})>~ 0$ and
$\lim_{n\to\infty}D(x_{2n},x_{2n+1})=0$ we get\\
\begin{align*}
\limsup_{m,n\to\infty}\Big(1-\frac{K^2}{\varphi\Big(D(x_{2n+1},x_{2m+2})\Big)}\Big)^{-1}=\infty.\\
\end{align*}
which implies,
$\limsup_{m,n\to\infty}\varphi\Big(D(x_{2n+1},x_{2m+2})\Big)={K^2}^+$
and so $\limsup_{m,n\to\infty}D(x_{2n+1},x_{2m+2})=~0$,
 which is a contradiction. Hence $\limsup_{m,n\to\infty}D(x_{n},x_{m})=0$.
Now, since $\lim_{m,n\to\infty}D(x_{n},x_{m})=~0$, so $\{x_n\}$ is
Cauchy,Since $(X,D,K)$ is a complete $b$-metric-like space, the
sequence $\{x_n\}$ is convergent to $z$. Hence
\begin{align*}
\lim_{n\to\infty}D(x_n, z)=\lim_{n\to\infty}D(x_{2n},
z)=\lim_{n\to\infty}D(x_{2n+1}, z)=D(z, z)=\lim_{n\to\infty}D(x_n,
x_m)=0.
\end{align*}
Since T, S are onto, there exists $w_1\in X$ and $w_2\in X$ such
that $Tw_1=z$ and $Sw_2=z$.\\
we prove that $w_1=w_2=z$. Suppose to the contrary that $z\neq w_1$ and
$z\neq w_2$ , then we have
\begin{align*}
D(x_{2n},z)=D(Tx_{2n+1},Sw_2)\geq
\varphi\Big(D(x_{2n+1},w_2)\Big)D(x_{2n+1},w_2).\\
D(x_{2n+1},z)=D(Sx_{2n+2},Tw_1)\geq
\varphi\Big(D(x_{2n+2},w_1)\Big)D(x_{2n+2},w_1).
\end{align*}
By taking limit as $n\to \infty$ in the above inequalities and
applying Lemma \ref{lem2} we have
 \begin{align*}
0=\lim_{n\to\infty}D(x_{2n},z)&\geq
\lim_{n\to\infty}\varphi\Big(D(x_{2n+1},w_2)\Big)\lim_{n\to\infty}D(x_{2n+1},w_2)\\&\geq
K^{-1} \lim_{n\to\infty}\varphi\Big(D(x_{2n+1},Z)\Big)D(Z,w_2),
\end{align*}
and
\begin{align*}
0=\lim_{n\to\infty}D(x_{2n+1},z)&\geq
\lim_{n\to\infty}\varphi\Big(D(x_{2n+2},w_1)\Big)\lim_{n\to\infty}D(x_{2n+2},w_1)\\&\geq
K^{-1} \lim_{n\to\infty}\varphi\Big(D(x_{2n+2},Z)\Big)D(Z,w_1),
\end{align*}
and hence $\lim_{n\to\infty}\varphi\Big(D(x_{2n+1},z)\Big)=0 $ and
$\lim_{n\to\infty}\varphi\Big(D(x_{2n},z)\Big)=0$ which is a
contradiction. Since $\varphi{t}>K^{2}$, for all $ t\in[0,\infty) $,
so $\lim_{n\to\infty}\varphi\Big(D(x_{2n},z)\Big)\geq K^2>0 $ and
$\lim_{n\to\infty}\varphi\Big(D(x_{2n+1},z)\Big)\geq K^2>0$.
Therefore
 $z=w_1=w_2$ that is $z=Tw_1=Tz=Sw_2=Sz$.
\end{proof}
We know that $b$-metric-like spaces are an extention of partial
metric, metric-like and $b$-metric spaces. Therefore we get the
following results.
\begin{corollary}
Let $(X,P)$ be a complete partial metric space.  Suppose $S,T:X\to
X$ are onto mappings,  such that  $P(Tx,Sy)\geq
\varphi(P(x,y))P(x,y)$ for all $x,y\in X$, then $T$ and $S$ have a
unique common fixed point.
\end{corollary}
\begin{corollary}
let $(X,\sigma)$ be a complete metric-like space. Suppose $S,T:X\to
X$ are onto mappings,  such that  $\sigma(Tx,Sy)\geq
\varphi(\sigma(x,y))\sigma(x,y)$ for all $x,y\in X$, then $T$and $S$
have a common fixed point.
\end{corollary}
\begin{corollary}
let $(X,d,K)$ be a complete  $b$-metric space, Assume that $S,T:X\to
X$ are onto  mappings,  such that  $D(Tx,Sy)\geq
\varphi\Big(d(x,y)\Big)d(x,y)$ for all $x,y\in X$, then $T$ and $S$
have a common fixed point.
\end{corollary}

\end{document}